\documentclass[11pt]{amsart}

\usepackage{tikz}
\usepackage{subfig}
\usepackage{amsmath,amssymb,cancel}
\usepackage[none]{hyphenat}
\usepackage{algorithm}
\usepackage{algpseudocode}
\usepackage{url}

\newcommand{\sk}{s_k}

\newtheorem{thm}{Theorem}[section]

\newtheorem{lem}[thm]{Lemma}

\newtheorem{conjecture}[thm]{Conjecture}

\numberwithin{equation}{section}
\numberwithin{algorithm}{section}

\begin{document}
\bibliographystyle{plain}

%Author and title info
\author{Sara Moore}
\address{Department of Physics \&\ Astronomy, University of Rochester, Rochester, New York 14627 }
\email{smoore55@u.rochester.edu}

 \author{Jonathan P. Sorenson}
 \address{Computer Science and Software Engineering Department,
   Butler University, Indianapolis, Indiana 46208}
\email{jsorenso@butler.edu}

\title{Explicit Bounds and Parallel Algorithms for Counting Multiply Gleeful Numbers}

\begin{abstract}
  Let $k\ge 1$ be an integer.
    A positive integer $n$ is $k$-\textit{gleeful} if $n$
    can be represented as the sum of $k$th powers of
    consecutive primes.
    For example, $35=2^3+3^3$ is a $3$-gleeful number,
  and $195=5^2+7^2+11^2$ is $2$-gleeful.
    In this paper, 
    we present some new results on $k$-gleeful numbers
    for $k>1$.
    
    First, we extend previous analytical work.
    For given values of $x$ and $k$, 
      we give explicit upper and lower bounds on the number of
      $k$-gleeful representations of integers $n\le x$.
      
    Second, we describe and analyze two new, efficient parallel algorithms,
    one theoretical and one practical, 
    to generate all $k$-gleeful representations up to a bound $x$.
    
    Third, we study integers that are
    \textit{multiply} gleeful, that is, integers with more than one representation as a sum of powers of
    consecutive primes, 
    including both the same or different values of $k$.
    We give a simple heuristic model for estimating
    the density of multiply-gleeful numbers,
    we present empirical data in support of our heuristics,
    and offer some new conjectures.
\end{abstract}

\date{\today}
\maketitle

\section{Introduction}\label{sec:intro}

Let $k\ge1$ be an integer.
We say a positive integer $n$ is $k$-\textit{gleeful} if $n$
  can be written as the sum of $k$th powers of consecutive
  primes.
For example, $35=2^3+3^3$ is a $3$-gleeful number,
  and $195=5^2+7^2+11^2$ is $2$-gleeful.
Let $f_k(n)$ denote the number of representations a positive
  integer $n$ has as a $k$-gleeful number,
  and let 
  $$ \sk(x) = \sum_{n=1}^x f_k(n),$$
  the total number of $k$-gleeful representations up to $x$.
In this paper, we address two questions about gleeful numbers
  and their representations:
  \begin{itemize}
      \item Can we give explicit upper and lower bounds on
      $\sk(x)$?
      \item What can we say about integers where $f_k(n)>1$?
  \end{itemize}
Before we state our results, we give some background on
  what is already known.

\subsection{Previous Work}

Moser \cite{Moser63} proved that $s_1(x) \sim x\log 2$.
He also posed several interesting questions on the behavior
  of $f_1(n)$.
See also \cite{UPINT}.
In this paper, we only look at $k$-gleeful numbers for integers $k>1$.

Tongsomporn, Wananiyakul, and Steuding \cite{TWS2022}
  proved that 
  $$s_2(x) < 10.9558 \frac{x^{2/3}}{(\log x)^{4/3}}.$$
They also computed a list of $2$-gleeful numbers up to 2000.

In \cite{OSS2024} it was proved that for every integer $k>1$,

\begin{equation}\label{eq:sopplower}
s_k(x) \ge \frac{(k+1)^2}{2} \cdot
 \frac{x^{2/(k+1)}}{(\log x)^{{2k/(k+1)}}}\cdot(1+o(1))
\end{equation} 
and for $c_k= (k^2/(k-1)) (k+1)^{1-1/k}$,
\begin{equation}\label{eq:soppupper}
s_k(x) \le c_k \cdot
\frac{x^{2/(k+1)}}{(\log x)^{{2k/(k+1)}}}(1+o(1)).
\end{equation}
Note that $(k+1)^2>c_k>(k+1)^2/2$.
They also gave two efficient sequential algorithms:
  one to enumerate $k$-gleeful representations
  and one to compute the exact value of $s_k(x)$,
  and they presented numerical data supporting their analytical results.

\subsection{New Results and Paper Outline}

In this paper, we continue the work from \cite{OSS2024}.

In \S\ref{sec:explicit}, for given values of
  $k$ and $x$, we give explicit upper and lower bounds
  on $s_k(x)$.

Our particular interest was to learn more
  about multiply-gleeful numbers or  \textit{duplicates}, 
  that is, integers $n$ with either $f_k(n)>1$
  or both $f_k(n)>0$ and $f_{k^\prime}(n)>0$ for $k\ne k^\prime$.
 
We found the enumeration algorithm from \cite{OSS2024} did not work
  well for finding duplicates, as it requires too much memory.
In \S\ref{sec:alg}, we describe and analyze two parallel
  algorithms for finding $k$-gleeful numbers -- one practical,
  the other theoretical.
The practical algorithm is based on
  a sequential routine that finds all $k$-gleeful numbers in a short interval.
The results from that interval are then sorted to detect
  values of $n$ with $f_k(n)>1$.
This parallelizes nicely by simply processing short intervals
  concurrently.
To detect duplicates with differing $k$ values,
  the algorithm is run twice on the same interval, once for each $k$ value, and again, the interval's results are sorted to detect the duplicates.
Our theoretical algorithm shows this problem is in $\mathcal{NC}$
  \cite{GHR}.

Then, in \S\ref{sec:dups}, we describe a heuristic model
  that predicts how many duplicates we expect to find up to $x$.
We then evaluated our model using data generated by our
  new parallel algorithm.
We state some conjectures consistent with these results.

Our code and data can be found at \url{https://github.com/sorenson64/sopp}.

\section{Explicit Bounds}\label{sec:explicit}

\newcommand{\skm}{{s_{k,m}}}
\newcommand{\Sk}{{s_k}}

%\begin{center}\textbf{Explicit Bounds for Sums of Consecutive Primes to the kth Power} \\ Sara Moore
%\end{center}

In this section, we prove the following theorems, 
which are explicit versions of Theorem 1 from \cite{OSS2024}.

Let $p_n$ denote the $n$th prime with $p_1=2$.
The number of primes up to $x$ is given by $\pi(x)$.
For fixed $x$ and $k\ge2$ an integer,
let $M:=M(x,k)$ be the maximum length of any representation of any
  integer up to $x$, so that
\begin{equation} \label{eq:x}
2^k + 3^k + \cdots + p_M^k \le x
  < 2^k + 3^k + \cdots + p_M^k+p_{M+1}^k.
\end{equation}
Observe that for a given $k$, 
any bound on $M$ implies a bound on $x$.
Let $M_0$ be a fixed integer, at least $6$.
Our results depend on $M(x,k)$ being larger than $M_0$.
Choosing $M_0$ to be larger gives better explicit constants.

%  of Theorem 1 and Lemma 2 from \cite{OSS2024}.
We define the following functions on $y$ and $k$.
We will be plugging $M_0$ in for $y$.
\begin{eqnarray*}
c_k& =& \left(\frac{k^2}{k-1}\right) \cdot (k+1)^{(k-1)/k} \quad
  \mbox{  from equation (\ref{eq:soppupper})},\\
    A(y)&:=& \frac{\log (y/2)}{\log y} , \\
    B(y,k)&:=& \frac{\log (y+1)}{\log y}+
  \frac{\log\log(y+1)^2}{\log y} \frac{k}{k+1}, \\
     C(y,k)&:=&
    \left(\frac{y}{y-1}\right)^{1/(k+1)} B(y,k)^{k/(k+1)}, \\
    D(y,k)&:=&
    \left(\frac{y}{y+3}\right)
 \left(\frac{\log(y/2)}{\log y + 2\log\log(y+2)}\right)^{k/(k+1)},\\
 E(y,k)&:=&1+\frac{1}{(k+1)A(y)-1},\\
F(y,k)&:=&\left(\frac{y+1}{y} \right)^{(k-1)/k} \cdot 4^{(k-1)/(k(k+1))}
  \cdot C(y,k)^{(k-1)/k}\cdot E(y,k), \\
  U(y,k) &:=& 1.25506 \cdot F(y,k), \\
  L(y,k) &:=& \left(\frac{y-1}{y} \right)\cdot D(y,k)^2.
\end{eqnarray*}
Note that $A(y),B(y,k),C(y,k),D(y,k)\rightarrow 1$ for large $y$.
$E(y,k)\rightarrow 1+1/k$ for large $y$.

%%%%%%%%%%%%%%%%%%%%%%%%%%%%%%%%%%%%%%%%%%%%%%%%%%%%%

\begin{thm}\label{thm:skxupper} For $k>1$ and $M\ge M_0\ge 6$, we have 

$$
s_k(x) \le  c_k \cdot
\frac{x^{2/(k+1)}}{(\log x)^{k/(k+1)}}
\cdot U(M_0,k).
$$
\end{thm}
For $k=2$, in \cite{TWS2022} they give $10.9558$.
Here, for $k=2$, we get the weaker bound $14.2423$ for $M\ge M_0=6$.
The results in \cite{OSS2024} give a constant
of $c_2=4\cdot\sqrt{3}\approx 6.928$, for large $x$

%$$S_k(x) \leq \p{\frac{3\p{2^{\frac{2-k}{k}}}k^3(k+1)^{\frac{3k-1}{k}}}{(k-1)^3}}\p{\frac{x^{\frac{2}{k+1}}}{(\log x) ^{\frac{2k}{k+1}}}}$$

\begin{thm}\label{thm:skxlower} For $k>1$ and $M\ge M_0\ge 6$, we have
$$ s_k(x) \ge \frac{(k+1)^2}{2} \cdot 
\frac{x^{2/(k+1)}}{(\log x)^{k/(k+1)}} \cdot L(M_0,k).$$
%$$ S_k(x) \geq \frac{15000(k+1)^2}{51005} \p{\frac{x^{\frac{2}{k+1}}}{\log^{\frac{2k}{k+1}} x}}$$

%which applies for  $M>5$
\end{thm}
See Table \ref{table:bounds}; we show two numbers for each combination of
$M_0$ and $k$: the lower bound constant followed by the upper bound constant.
These constants include everything except the main term
$x^{2/(k+1)}/(\log x)^{k/(k+1)}$.

\begin{table}[h]
    \centering
\begin{tabular}{l|c|c|c|c|}
%\begin{tabular}{l|r|r|r|r|r|}
$M_0=$ & $6$ & $100$ & $10000$ & $1000000$  \\ \hline
$k=2$  
 & 0.391504, 14.2423
 & 1.71182, 12.1097
 & 2.39745, 11.6778
 & 2.7343, 11.5116
 \\ 
$k=3$  
 & 0.580731, 23.4232
 & 2.72032, 18.7705
 & 3.93987, 17.7299
 & 4.5675, 17.3147
 \\ 
$k=5$  
 & 1.09023, 63.156
 & 5.47127, 48.0799
 & 8.19445, 44.4013
 & 9.6564, 42.8989
 \\ 
$k=10$  
 & 3.10821, 249.625
 & 16.6068, 182.224
 & 25.6426, 164.599
 & 30.6698, 157.311
 \\ 
$k=20$  
 & 10.3113, 1009.68
 & 57.0995, 720.629
 & 89.7176, 642.315
 & 108.222, 609.797
 \\ \hline
 \end{tabular}
    \caption{Constants for lower and upper bounds on $\Sk(x)$}
    \label{table:bounds}
\end{table}
%%%%%%%%%%%%%%%%%%%%%%%%%%%%%%%%%%%%%%%%%%%%%%%%%%%%%

\subsection{Setup}

Let us define $f_{k,m}(n)$ to be the number of representations
  of $n$ as a sum of exactly $m$ $k$th powers of consecutive primes.
Observe that $f_{k,m}(n)=0$ or $1$ and that
$$
f_{k}(n)=\sum_{m=1}^M f_{k,m}(n).
$$
Let $\skm(x)$ be the number of positive integers $n$ such
  that
$$ p_n^k+p_{n+1}^k+\cdots+p_{n+m-1}^k\le x.$$
Then
$$
\skm(x)=\sum_{n=1}^x f_{k,m}(n)
$$
so that
$$ \Sk(x)=\sum_{n=1}^x f_{k}(n)=\sum_{n=1}^x\sum_{m=1}^M f_{k,m}(n)=
\sum_{m=1}^{M} \skm(x).$$
This counts the number of representations of $k$-gleeful
numbers $\le x$.

%%%%%%%%%%%%%%%%%%%%%%%%%%%%%%%%%%%%%%%%%%%%%%%%%%%%%

We will make use of the following results due to Rosser and Schonfeld \cite{RS62} and Rosser \cite{Rosser1939}.
\begin{eqnarray}
  n\log n  <  p_n &&\quad\mbox{ for } n\ge 1
  \label{eq:pnlow}\\
  p_n < n\log n + 2n\log\log n &&\quad\mbox{ for } n\ge 3
  \label{eq:pnhigh}\\
  p_n < 2n\log n && \quad\mbox{ for }n\ge3
  \label{eq:pnhighlame}\\
  \pi(x)<1.25506\frac{x}{\log x}
  && \quad\mbox{ for }x\ge 2
  \label{eq:pixhigh}
\end{eqnarray}
%%%%%%%%%%%%%%%%%%%%%%%%%%%%%%%%%%%%%%%%%%%%%%%%%%%%%
We will need upper and lower bounds on $M$.

\begin{lem}\label{lem:mhigh}
If $M\ge M_0 \ge 6$ then
$$
% M < 8.856^{1/(k+1)}\cdot(k+1)\frac{ x^{1/(k+1)}}{(\log x)^{k/(k+1)}}
% < 2.069 \cdot(k+1)\frac{ x^{1/(k+1)}}{(\log x)^{k/(k+1)}}
M < 4^{1/(k+1)} \cdot \frac{(k+1) x^{1/(k+1)}}{(\log x)^{k/(k+1)}}
\cdot C(M_0,k)
.
$$
\end{lem}

\begin{lem}\label{lem:mlow}
If $M\ge M_0\ge 6$ then
$$
M \ge (k+1)\cdot\frac{ x^{1/(k+1)}}{(\log x)^{k/(k+1)}}
\cdot D(M_0,k).
$$
%$M$ is bounded below by
%
%$$M > \p{\frac{100}{101}}\p{\frac{3}{4}}^{\frac{k}{k+1}}\p{k+1}\frac{x^{\frac{1}{k+1}}}{\log^{\frac{k}{k+1}}x}$$
\end{lem} 
For large $M$, from \cite{OSS2024} we expect 
$M\sim (k+1) \cdot x^{1/(k+1)}/(\log x)^{k/(k+1)}$.
See Table \ref{table:1} for some exact values of $M(x,k)$.

%%%%%%%%%%%%%%%%%%%%%%%%%%%%%%%%%%%%%%%%%%%%%%%%%%%%%
\begin{table}[ht]
\centering
 \begin{tabular}{|c |c |c| c| c| c|} 
 \hline
 x & $M(x,2)$ & $M(x,3)$ & $M(x,5)$ & $M(x,10)$ & $M(x,20)$\\ [0.5ex] 
 \hline
 $10^{3}$  & 7 & 4 & 2 & 0 & 0 \\
 $10^{4}$  & 14 & 7 & 3 & 1 & 0 \\
 $10^{5}$  & 28 & 11 & 4 & 2 & 0 \\
 $10^{6}$  & 54 & 18 & 6 & 2 & 0 \\
 $10^{7}$  & 105 & 29 & 8 & 3 & 1 \\
 $10^{8}$  & 207 & 47 & 11 & 3 & 1 \\
 $10^{9}$  & 411 & 77 & 15 & 4 & 1 \\
 $10^{10}$ & 822 & 126 & 21 & 4 & 2 \\
 $10^{11}$ & 1656 & 209 & 30 & 5 & 2 \\
 $10^{12}$ & 3356 & 348 & 40 & 6 & 2 \\
 $10^{13}$ & 6834 & 581 & 55 & 8 & 2 \\
 $10^{14}$ & 13975 & 974 & 76 & 9 & 3 \\
 $10^{15}$ & 28682 & 1640 & 106 & 10 & 3 \\
 $10^{16}$ & 59066 & 2771 & 148 & 12 & 3 \\
 $10^{17}$ & 121987 & 4695 & 206 & 15 & 4 \\
 $10^{18}$ & 252574 & 7977 & 288 & 17 & 4 \\
 $10^{19}$ & 524136 & 13589 & 403 & 20 & 4 \\
 $10^{20}$ & 1089888 & 23201 & 566 & 24 & 4 \\
 \hline
 \end{tabular}
 \caption{Exact values of $M(x,k)$ for various $x,k$}
\label{table:1}
\end{table}

%%%%%%%%%%%%%%%%%%%%%%%%%%%%%%%%%%%%%%%%%%%%%%%%%%%%%
As stepping stones, 
the proofs of these two lemmas utilize easier-to-prove
upper and lower
bounds on $\log M$ in terms of $x$ and $k$.

\begin{lem}\label{lemma:logm}
If $M\ge M_0\ge 6$ then we have
\begin{eqnarray*}
\log M &<& \frac{\log x}{k+1} \cdot \frac{1}{A(M_0)} \quad \mbox{and}\\
\log M &\ge& \frac{\log x}{k+1}\cdot \frac{1}{B(M_0,k)} .
% \\
%  \frac{1}{\log M}&>&\frac{k+1}{\log x}\cdot
% \frac{\log (M_0/2)}{\log M_0} \quad\mbox{and} \\
% \frac{1}{\log M}&\le& \frac{k+1}{\log x}\cdot
%   \left(\frac{\log (M_0+1)}{\log M_0}+
%   \frac{\log\log(M_0+1)^2}{\log M_0} \frac{k}{k+1}
%   \right) .
\end{eqnarray*}
\end{lem}
As we let $M_0$ get larger, we get the expected
$(k+1)\log M\sim \log x$.
\begin{proof}
Since $M\ge 6$, we can bound the sum on the left of (\ref{eq:x}) from below by $(M/2) p_{M/2}^k$.
Taking $(M/2)\log (M/2)<p_{M/2}$ from (\ref{eq:pnlow}), we have
$$
(M/2)^{k+1} ( \log (M/2) )^k < x.
$$
Because $M\ge 6 > 2e$, we can drop the
  $(\log (M/2))^k$ term since it exceeds 1.
Taking logarithms of both sides then gives 
\begin{eqnarray*}
\log x > (k+1)\log (M/2) 
&=& (k+1)(\log M)\left(1 - \frac{\log 2}{\log M}\right)\\
&\ge&(k+1)(\log M)\left(1 - \frac{\log 2}{\log M_0}\right).
\end{eqnarray*}

The sum on the right of (\ref{eq:x}) is bounded above by $(M+1) p_{M+1}^k$.
Using (\ref{eq:pnhighlame}) gives
$$
\log x < (k+1)\log(M+1)+ k\log (2\log (M+1)).
$$
With $M\ge M_0$, we know that 
\begin{equation} \label{eq:logm1}
\log(M+1)= \frac{\log (M+1)}{\log M}\log M 
  < \frac{\log (M_0+1)}{\log M_0}\log M.
\end{equation}
We also know that when $M\ge M_0\ge 6$,
$$
\frac{\log (2\log (M+1))}{\log (M+1)}
$$
is maximized when $M=M_0$.
This then gives
\begin{eqnarray*}
\frac{\log x}{k+1} &\le&  (\log M)
  \left(1+ \frac{\log\log(M_0+1)^2}{\log (M_0+1)} \frac{k}{k+1}  \right)\frac{\log (M_0+1)}{\log M_0} \\
  &=&(\log M)\left(\frac{\log (M_0+1)}{\log M_0}+
  \frac{\log\log(M_0+1)^2}{\log M_0} \frac{k}{k+1}
  \right).
\end{eqnarray*}
Plugging in the definitions of $A(M_0),B(M_0,k)$ completes the proof.
%Observing that $k/(k+1)<1$ completes the proof.
\end{proof}

%%%%%%%%%%%%%%%%%%%%%%%%%%%%%%%%%%%%%%%%%%%%%%%%%%%%%
\subsection{Proof of Lemma \ref{lem:mhigh}}

Next, we bound $x$ from below in terms of $M$ and $k$.
Using (\ref{eq:x}) and (\ref{eq:pnlow}), We have
\begin{eqnarray*}
x &\ge& p_1^k + p_2^k + \cdots + p_M^k \\
  & > & \sum_{n=1}^{M}(n \log n)^k 
  \quad \ge\quad  \sum_{n=M^{1-1/k}}^{M}(n \log n)^k \\
  & \ge & (\log M^{1-1/k})^k \sum_{n=M^{1-1/k}}^{M}n^k \\
  & = & (1-1/k)^k (\log M)^k \sum_{n=M^{1-1/k}}^{M}n^k. \\
\end{eqnarray*}
We can bound the sum on $n^k$ with an integral to get
\begin{eqnarray*}
\sum_{n=M^{1-1/k}}^{M}n^k
 &\ge& 
  \int_{M^{1-1/k}}^{M} t^k dt \quad = \quad
   \frac{M^{k+1}-{M^{(k+1)(1-1/k)}}}{k+1} \\
  &=& \frac{M^{k+1}-{M^{k-1/k}}}{k+1}
  \quad = \quad
  \frac{M^{k+1}}{k+1} \left(1 -\frac{1}{M^{1+1/k}}\right) \\
  &>& \frac{M^{k+1}}{k+1}\cdot \frac{M_0-1}{M_0}
\end{eqnarray*}
if $M\ge M_0$.
We also have $(1-1/k)^k \ge 1/4$ for $k\ge 2$.
Pulling this together gives
$$
x \ge \frac{M^{k+1} (\log M)^k}{4(k+1)}\cdot \frac{M_0-1}{M_0},
$$
which is valid when $M\ge M_0$.
This directly gives
$$ M^{k+1}\le \frac{4(k+1) x}{(\log M)^k}\cdot \frac{M_0}{M_0-1} .$$
Applying Lemma \ref{lemma:logm} then gives
$$
M^{k+1}< \frac{4\cdot(k+1)^{k+1}  x}{(\log x)^k}
  \cdot \frac{M_0}{M_0-1} B(M_0,k)^k
$$
and
$$
M < 4^{1/(k+1)}\frac{(k+1) x^{1/(k+1)}}{(\log x)^{k/(k+1)}}
\cdot \left(\frac{M_0}{M_0-1}\right)^{1/(k+1)} B(M_0,k)^{k/(k+1)}.
$$
This completes the proof of Lemma \ref{lem:mhigh}.

%Because $p_n > n\log n$ (Rosser \& Schoenfeld, 1962), 
%
%$$x \geq p_1^k + p_2^k + \cdots + p_M^k > \sum_{n=2}^{M}\p{n \log n}^k \geq 
%\sum_{n=2}^{\sqrt{M}} n^k + \sum_{n=\sqrt{M}}^{M}n^k \log ^kn \geq
%\sum_{n=2}^{\sqrt{M}} n^k + \log^k\sqrt{M}
% \sum_{n=\sqrt{M}}^{M}n^k$$
% $$\geq \frac{1}{2}\log^k M \sum_{n=\sqrt{M}}^{M}n^k \geq \int_{\sqrt{M}-1}^{M}n^k dn = \frac{\log ^kn}{2(k+1)}\p{M^{k+1}-(\sqrt{M}-1)^{k+1}} \geq 
% \frac{\log ^kn}{2(k+1)}\p{M^{k+1}-(\sqrt{M})^{k+1}}$$
%
% $$> \frac{\log ^kM}{2(k+1)}\frac{M^{k+1}}{2}$$

%Then because $$x > \frac{M^{k+1}\log ^kM}{4(k+1)} > M^{k+1}$$ for $M > 32$, then $\log x > (k+1)\log M$
%or $$\log^k M < \frac{\log^k x}{(k+1)^k} $$

%Applying this and solving this inequality for $M$ results in $$M < \frac{4^{\frac{1}{k+1}}(k+1)x^{\frac{1}{k+1}}}{\log ^{\frac{k}{k+1}}x}$$

%%%%%%%%%%%%%%%%%%%%%%%%%%%%%%%%%%%%%%%%%%%%%%%%%%%%%%%

\subsection{Proof of Lemma \ref{lem:mlow}}

Again, starting from (\ref{eq:x}),
observing that $2^k+3^k<5^k<p_{M+2}^k$, and applying (\ref{eq:pnhigh}), we have
\begin{eqnarray*}
    x &\le& p_1^k+\cdots+p_{M+1}^k 
    \le p_3^k+\cdots+p_{M+2}^k \\
    % &\le& 2\left( 
    % p_{\lfloor (M+1)/2 \rfloor}^k+\cdots+p_{M+1}^k
    % \right) \quad=\quad 2 \sum_{n=\lfloor(M+1)/2 \rfloor}^{M+1} p_n^k \\
    &\le&  \sum_{n=3 \rfloor}^{M+2}
      n^k ( \log n + 2\log\log n)^k \\
    &\le& (\log (M+2) + 2\log\log (M+2))^k
      \sum_{n=3}^{M+2} n^k \\
    &\le& (\log (M+2) + 2\log\log (M+2))^k 
      \frac{(M+3)^{k+1}}{k+1} \\
    &\le&  (\log M)^k 
    \left( \frac{\log(M_0+2)}{\log M_0}+\frac{2\log\log(M_0+2)}{\log M_0}\right)^k
    \frac{(M+3)^{k+1}}{k+1}.
\end{eqnarray*}
Here we bounded the sum $\sum_{n=3}^{M+2} n^k$
with an integral to get $\frac{(M+3)^{k+1}}{k+1}$.
% Observe that $2\log \log t / \log t \le 2/e$ and
% that $(M+2)/M\le 4/3$, and  then by
Using (\ref{eq:logm1}), we have
$$
x \le ( \log M )^k \frac{((M_0+3)/M_0)\cdot M)^{k+1}}{k+1}
\left( \frac{\log(M_0+2)}{\log M_0}+\frac{2\log\log(M_0+2)}{\log M_0}\right)^k,
$$
or
$$
(k+1) \frac{x}{(\log M)^k} 
\frac{(M_0/(M_0+3))^{k+1}}{ \left( \frac{\log(M_0+2)}{\log M_0}+\frac{2\log\log(M_0+2)}{\log M_0}\right)^k}
\le M^{k+1}.
$$
We apply Lemma \ref{lemma:logm} and simplify a bit to obtain
$$
 \frac{(k+1) x^{1/(k+1)}}{(\log x)^{k/(k+1)}} \cdot
 \left(\frac{M_0}{M_0+3}\right)
 \left(\frac{\log(M_0/2)}{\log( M_0+2) + 2\log\log(M_0+2)}\right)^{k/(k+1)}
 \le M.
$$
This completes the proof.

%Noting that $$p_n < n(\log n + \log \log n)$$ 
%
%(Rosser \& Schoenfeld, 1962),
%
%$$p_1^k + p_2^k + \cdots + p_M^k < \sum_{m=2}^M (n \log n + n\log \log n)^k = \sum_{m=2}^M \p{n \log n}^k \p{1 + \frac{\log \log n}{\log n}}^k$$
%
%$$ < \sum_{m=2}^M \p{\frac{4}{3}}^k n^k \log^k n \leq \p{\frac{4}{3}}^k \log^k M \sum_{m=2}^M n^k < \p{\frac{4}{3}}^k \log^k M \int_2^{M+1} n^k dk$$
%
%
%$$= \p{\frac{4}{3}}^k \log^k M \p{\frac{(M+1)^{k+1}}{k+1} - \frac{2^{k+1}}{k+1}}
%< \p{\frac{4}{3}}^k \p{\frac{1}{k+1}} \log^k M \p{(M+1)^{k+1}}$$ 
%$$\leq \p{\frac{4}{3}}^k \p{\frac{1}{k+1}} \log^k M \p{(\frac{101M}{100})^{k+1}} = \p{\frac{4}{3}}^k \p{\frac{101}{100}}^{k+1}\p{\frac{1}{k+1}} \log^k M \p{M^{k+1}}$$ 
%
%Also noting that $$\p{\frac{4}{3}}^k \p{\frac{101}{100}}^{k+1}\p{\frac{1}{k+1}} \log^k M \p{M^{k+1}} < M^{k+2}$$ for $M > 688$, we can solve this inequality for M, resulting in 
%
%$$M > \p{\frac{100}{101}}\p{\frac{3}{4}}^{\frac{k}{k+1}}\p{k+1}^{\frac{1}{k+1}}\p{k+2}^{\frac{k}{k+1}}\frac{x^{\frac{1}{k+1}}}{\log^{\frac{k}{k+1}}x}
%> \p{\frac{100}{101}}\p{\frac{3}{4}}^{\frac{k}{k+1}}\p{k+1}\frac{x^{\frac{1}{k+1}}}{\log^{\frac{k}{k+1}}x}
%$$

%%%%%%%%%%%%%%%%%%%%%%%%%%%%%%%%%%%%%%%%%

\subsection{Proof of Theorem \ref{thm:skxupper}}
For positive integers $n,m$, we refer to any sum of the form
$p_n^k + p_{n+1}^k + \cdots + p_{n-1+m}^k$
as a \textit{chain} of length $m$.
Recall that $$\skm(x) = \#\{n:p_n^k + p_{n+1}^k + \cdots + p_{n-1+m}^k \leq x\},$$
the number of $k$-gleeful representations of integers $\le x$ of length $m$,
or the number of chains of length $m$ whose sums are bounded by $x$.

For the moment, fix a chain length $m$.  
Choose $n$, the starting point of the chain, as large as possible so that we have
$$m\cdot p_n^k \leq p_n^k + p_{n+1}^k + \cdots + p_{n-1+m}^k \leq x \leq  p_{n+1}^k + p_{n+2}^k + \cdots + p_{n+m}^k.$$
This gives then
\begin{eqnarray*}
    mp_n^k &\le& x, \\
    p_n &\le& (x/m)^{1/k}, \\
    n &\le& \pi( (x/m)^{1/k} ) .
\end{eqnarray*}
Observe that $\skm(x) = n$ here.  (See Lemma 1 from \cite{OSS2024}.)
Thus,
$$
\Sk(x) = \sum_{m=1}^M \skm(x) \le \sum_{m=1}^M \pi( (x/m)^{1/k} ).
$$
Observe that $(x/m)^{1/k}\ge ( p_M^k/M )^{1/k}>2$ for $M\ge 6$.
From (\ref{eq:pixhigh}) we have %\cite[equation (3.6)]{RS62} we have
$$
\pi(t) \le 1.25506\cdot t/\log t
$$
when $t\ge 2$.
This gives
\begin{eqnarray*}
    s_k(x) &\le& \sum_{m=1}^M \pi( (x/m)^{1/k} )\\
    &\le& 1.25506 \cdot \sum_{m=1}^M \frac{k (x/m)^{1/k}}{\log(x/m)} \\
    &\le& 1.25506 \cdot \frac{k x^{1/k}}{\log(x/M)} \sum_{m=1}^M \frac{1}{m^{1/k}}.
\end{eqnarray*}
Focusing first on the logarithm in the denominator, we have
% \begin{eqnarray*}
%     \frac{1}{\log(x/M)}& =& \frac{1}{\log x - \log M} \\
%     &\le& \frac{1}{\log x - (1.632/(k+1)) \log x} \quad = \quad
%     \frac{k+1}{(k-0.632)\log x}
% \end{eqnarray*}
\begin{eqnarray*}
  \log(x/M) &\ge& \log x - \frac{\log x}{(k+1)A(M_0)} \quad\mbox{ or }\\
    \frac{1}{\log(x/M)} &\le& 
      \frac{1}{\log x} \left(1+\frac{1}{(k+1)A(M_0)-1}\right) \\
      &=& \frac{1}{\log x} \cdot E(M_0,k)
\end{eqnarray*}
using Lemma \ref{lemma:logm}.
Next, we estimate the sum:
\begin{eqnarray*}
    \sum_{m=1}^M \frac{1}{m^{1/k}} &\le& \int_1^{M+1} t^{-1/k} dt \\
      & \le & \frac{(M+1)^{1-1/k}}{1-1/k}.
\end{eqnarray*}
Pulling this together, we have
\begin{eqnarray*}
s_k(x) &\le& 1.25506 \cdot E(M_0,k)\cdot \frac{kx^{1/k}}{\log x} \frac{(M+1)^{1-1/k}}{1-1/k} \\
 %& = & 1.25506 \cdot \left(\frac{k^2(k+1)}{(k-1)(k-0.632)}\right) \frac{x^{1/k}}{\log x} \cdot (M+1)^{(k-1)/k} \\
 &\le& 1.25506 \cdot E(M_0,k)
 \cdot \left(\frac{M_0+1}{M_0} \right)^{(k-1)/k}
 \left(\frac{k^2}{k-1}\right) 
 \frac{x^{1/k}}{\log x} \cdot M^{(k-1)/k},
\end{eqnarray*}
since $M\ge M_0$.
Next, we plug in our upper bound on $M$ from Lemma \ref{lem:mhigh}.
We have
\begin{eqnarray*}
    M^{(k-1)/k} &<& 4^{(k-1)/(k(k+1))}\cdot (k+1)^{(k-1)/k}\cdot 
    C(M_0,k)^{(k-1)/k}\cdot
    \frac{x^{(k-1)/(k(k+1))}}{(\log x)^{(k-1)/(k+1)}}.
\end{eqnarray*}
Plugging this in, we obtain
\begin{eqnarray*}
s_k(x) &\le& 1.25506\cdot c_k   \cdot
 \frac{x^{2/(k+1)}}{(\log x)^{2k/(k+1)}}\cdot F(M_0,y),
\end{eqnarray*}
and the result follows.

\subsection{Proof of Theorem \ref{thm:skxlower} }

Any subsequence sum of the maximum chain length $M$ represents a $k$-gleeful
number $\le x$.  Thus,
 the number of $i$, $j$ pairs such that $1\leq i \leq j \leq M$, 
or $\binom{M}{2}$, is a lower bound for $s_k(x)$.  Thus,
$$
s_k(x) \ge \binom{M}{2} = \frac{M(M-1)}{2} \ge \frac{M_0-1}{2M_0}\cdot M^2
$$
since we are assuming $M\ge M_0$.  
Simply apply Lemma \ref{lem:mlow} and
a bit of algebra and the result follows.

% Thus, by applying Theorem 2, we compute a lower bound of 

% $$S_k(x) \geq \binom{M}{2} = \frac{M(M-1)}{2} \geq \frac{2M^2}{5} \geq \frac{2}{5} \p{\p{\frac{100}{101}}\p{\frac{3}{4}}^{\frac{k}{k+1}}\p{k+1}\frac{x^{\frac{1}{k+1}}}{\log^{\frac{k}{k+1}}x}}^2 = \frac{10000\p{3^{\frac{2k}{k+1}}}(k+1)^2}{51005\p{2^{\frac{k-1}{k+1}}}} \p{\frac{x^{\frac{2}{k+1}}}{\log^{\frac{2k}{k+1}} x}}$$
% $$> \frac{15000(k+1)^2}{51005} \p{\frac{x^{\frac{2}{k+1}}}{\log^{\frac{2k}{k+1}} x}}$$

% ALL FIGURES AND TABLES HAVE YET TO BE UPDATED
% \begin{figure}[ht!]
% \centering
% \includegraphics[scale=0.6]{Figures/k=3_ul.png}
% \caption{
% \label{fig:k=3}}
% \end{figure}

% \begin{figure}[ht!]
% \centering
% \includegraphics[scale=0.6]{Figures/k=5_ul.png}
% \caption{
% \label{fig:k=5}}
% \end{figure}

% \begin{figure}[ht!]
% \centering
% \includegraphics[scale=0.6]{Figures/k=10_ul.png}
% \caption{
% \label{fig:k=10}}
% \end{figure}

\section{Two Parallel Algorithms}\label{sec:alg}

We begin with a straightforward adaptation of the
enumeration algorithm from \cite{OSS2024} to work on an
interval.

\subsection{An Algorithm to Enumerate Representations on an Interval}

    Here we describe an algorithm that generates all integers $n$
    with $f_k(n)>0$, where $x_1\le n < x_2$ for inputs $k,x_1,x_2$.
    We obtain the original algorithm from \cite{OSS2024} by
    setting $(x_1,x_2)=(1,x)$.

    Let $x$ be the largest value of $x_2$ we plan to use in any application of this algorithm.  
    As a preprocessing step, we find all primes up to $x^{1/k}$ and compute the prefix array $r[]$, where $r[0]=0$ and $r[j]=r[j-1]+p_j^k$ where $p_j$ is the $j$th prime with $p_1=2$.

For a particular value of $n$ with $f_k(n)>0$,
we write $n=r[t]-r[b]$, a difference of prefix sum values,
which gives its representation as
$p_{b+1}^k+\cdots +p_t^k$.
The trick is to generate exactly the correct values of $b$ and $t$ 
to ensure that $x_1\le n < x_2$.
The outer loop iterates through all possible $b$ values,
  and the inner loop iterates through the correct $t$ values.
Let $t_s$ indicate the smallest $t$ for a given $b$.
Observe that as $b$ increases, $t_s$ is non-decreasing 
and $t_s> b$.  
See Algorithm \ref{alg:seq}.
\begin{algorithm}[ht]
\caption{Enumerate integers $n$ with $x_1\le n<x_2$
  and $f_k(n)>0$}\label{alg:seq}
\begin{algorithmic}
\Require{Integers $k>1$, $x_1<x_2$, a list of all primes up to $x_2^{1/k}$, and the prefix sum array $r[\,]$}
\State{$t_s \gets 1$}
 \State $\ell\gets\pi(x_2^{1/k})$
\For {$b\gets 0$ to $\ell$ }
      \While {$t_s\le \ell$ and $t_s\le b$} 
        \State $t_s\gets t_s+1$
      \EndWhile
      \While{ $t_s\le \ell$ and $r[t_s]-r[b]<x_1$} 
       \State $t_s\gets t_s+1$
      \EndWhile
      \For{ $t\gets t_s$ to $\ell$} 
        \State $n\gets r[t]-r[b]$
        \If{$x_1\le n < x_2$}
          \State output $(n,p_{b+1})$
        \ElsIf{$n\ge x_2$}
          \State {break the inner for-loop (on $t$)}
        \EndIf
      \EndFor 
    \EndFor 
\end{algorithmic}
\end{algorithm}
% \begin{quote}
%     \begin{tabbing}MM\=MM\=MM\=\kill
%     $t_s:=1$; \\
%     $\ell:=\lfloor\pi(x_2^{1/k})\rfloor$; \\
%     for $b=0$ to $\ell$ do:\+\\
%       while $t_s\le \ell$ and $t_s\le b$ do\+\\
%         $t_s:=t_s+1;$\-\\
%       endwhile; \\
%       while $t_s\le \ell$ and $r[t_s]-r[b]<x_1$ do\+\\
%         $t_s:=t_s+1;$\-\\
%       endwhile; \\
%       for $t:=t_s$ to $\ell$ do\+\\
%         $n:=r[t]-r[b]$; \\
%         if $x_1\le n < x_2$ then \+\\
%           output $(n,p_{b+1})$; \-\\
%         else if $n\ge x_2$ then \+\\
%           break the inner loop (on $t$); \-\\
%         endif; \-\\
%       endfor; \-\\
%     endfor; \\
%     \end{tabbing}
% \end{quote}

% Each $k$-gleeful number is represented as a difference
% from the prefix sum, $f[t]-f[b]$.
% Here $0\le b \le t \le \ell$.
% The idea is to generate all such differences in $[x_1,x_2]$
% while minimizing extra work.
% To do this, the outer loop runs through all possible values
% for $b$.
% $t_s$, then, is the smallest $t\ge b$ such that
% $f[t_s]-f[b]\ge x_1$.
% The first thing we do in the loop over $b$ values is to
% increment $t_S$ to meet this invariant.
% Observe that the total time spent incrementing $t_s$ is
% bounded by $\ell$.
% With $t_s$ set, we then run through all $t$ values, starting
% with $t_s$, generating all differences in our interval.
% The time spent on this part is at most proportional to
% $S_k(x_2)-S_k(x_1)$
% summed over all outer loop iterations.

The running time of this algorithm is bounded by a constant
times the number of times through the while loops and the
inner for loop.
Observe that the while loops increment $t_s$, which is bounded
  by $\ell = \pi(x^{1/k})$,
  so the total number of while-loop iterations is bounded by $\ell$.
The number of times we iterate through the inner for-loop is bounded by
a constant times the number of times the output ($n,p_{b+1})$
statement executes, which in turn is
$\Sk(x_2)-\Sk(x_1)$.

We have proven the following:
\begin{thm}
    Given integers $k>1$ and $x_1<x_2\le x$,
    Algorithm \ref{alg:seq} will list all integers $n$
    with $f_k(n)>0$ and $x_1\le n < x_2$.
    The number of arithmetic operations used by the algorithm
    is at most
    $O( x^{1/k}/\log\log x+(\Sk(x_2)-\Sk(x_1)))$.

    In addition, for every representation of $n$ as a 
    $k$-gleeful number, the first prime in that representation
    is also given.
\end{thm}
We have a few comments:
\begin{itemize}
    
\item
The understanding is that $x$ is an upper bound on the application
of the algorithm, and that it may be used on multiple intervals
$[x_1,x_2)$ with $x_1<x_2\le x$
We are also assuming here that $k$ is fixed, but $x$ (and $x_1,x_2$)
are large.

\item
The $x^{1/k}/\log\log x$ term is the time to compute
the list of primes up to $x^{1/k}$ using, say,
the Atkin-Bernstein algorithm \cite{AB2004}.
If the primes are already available, this term changes to
$\pi(x_2^{1/k})=O(kx^{1/k}/\log x)$.

\item 
In practice, we make the interval length $x_2-x_1$
large enough so that we expect $\Sk(x_2)-\Sk(x_1)\gg x^{1/k}$,
thereby ensuring that the cost of
managing the list of primes becomes negligible.

\item 
We obtain a
practical parallel algorithm by dividing the range $(1,x)$
into equal-sized intervals of length $\Delta=x_2-x_1$.
We than assign one processor to each interval.
Thus, $x/\Delta$ processors can run in parallel with little
communication overhead, and the list of primes and the prefix sum
array $r[\,]$ can be shared.
\item 
When searching for duplicates, which are integers $n$ with either
  $f_k(n)>1$ or both $f_k(n)>0$ and $f_{k^\prime}(n)>0$ for $k\ne k^\prime$,
  the interval of size $\Delta$ should be sorted to look for matches.
  Note that a list of integers in a limited range can be sorted
  in linear time using a radix or bucket-style sort.

  In practice, we found quicksort \cite{Hoare62} was good enough.
  See \cite[\S5.2.5]{Knuthv3}.
\end{itemize}

\subsection{A Theoretical Parallel Algorithm}

We describe the steps and analyze the algorithm as we go.
We assume an EREW PRAM parallel model with arithmetic operations
on integers with $O(\log x)$ bits taking constant time.
Note that $\pi(x^{1/k})=O(kx^{1/k}/\log x)$ by the prime number theorem.
As above, we use $\ell$ for the number of such primes.
\begin{enumerate}
    \item 
    To find the primes up to $x^{1/k}$, we use the algorithm from \cite{SP94}.
    This takes $O((1/k)\log x)$ time and 
    $O(kx^{1/k}/(\log x \log\log x))$
    processors.
    \item
    To compute the $k$th powers of all the primes,
    we use a sequential binary exponentiation algorithm that takes
    $O(\log k)=O(\log\log x)$ time,
    since we can assume $k=O(\log x)$ here.
    We apply this to all primes in parallel, 
    taking $\ell=O(kx^{1/k}/\log x)$ processors.
    \item 
    The prefix sum array $r[\,]$ can be computed 
    in $O(\log \ell)$ time using $O(\ell/\log \ell)$ processors, or $O((1/k)\log x)$ time and $O(k^2x^{1/k}/(\log x)^2)$ processors.
    \item 
    To start, we assign one processor to each $b$ value
    from $0$ to $\ell$.
    Then, for each $b$ in parallel, we perform a binary search
    on the $r[\,]$ array to find the correct start and stop
    $t$ values, $(t_1(b), t_2(b))$ so that for every $t$ 
    with $t_1(b)\le t\le t_2(b)$, we have
    $0<r[t]-r[b]<x$.
    This takes $O(\log \ell)$ time using $O(\ell)$ processors
    since this is how many $b$ values there are.
    \item 
    For each $b$, we allocate $(t_2(b)-t_1(b))/\log \ell$ additional
    processors.
    This is a total of $O(\ell+\Sk(x)/\log \ell)$ processors overall.
    \item 
    For every $b$, in parallel we compute $n=r[t]-r[b]$ for every
    $t_1(b)\le t \le t_2(b)$
    and output $(p_{b+1},n)$.
    This uses the processors allocated in the previous step.
    Each processor may have to do up to $O(\log \ell)$ such
    computations, but they take constant time each.
    This takes $O(\log \ell)$ time using
    $O( \ell + \Sk(x)/\log \ell)$ processors. 

\end{enumerate}
We have proven the following.
\begin{thm}
    There is an EREW PRAM algorithm to find all integers $n\le x$
    with $f_k(n)>0$ that uses at most $O(\log\ell)$ time
    and $O(\ell+s_k(x)/\log\ell)$ processors,
    where $\ell:=\pi(x^{1/k})=O(kx^{1/k}/\log x)$.
\end{thm}
Note that this algorithm is work-optimal, and proves that
computing $\sk(x)$ is in the complexity class $\mathcal NC$.

\section{Duplicates}\label{sec:dups}

In this section, we examine, heuristically, the distribution of duplicates, which come in two varieties:
  \begin{enumerate}
      \item Integers $n$ with $f_k(n)>1$, or
      \item Integers $n$ with both $f_k(n)>0$ and $f_{k^\prime}(n)>0$ for $k\ne k^\prime$.
  \end{enumerate}
WLOG, in the second type we shall henceforth assume $k<k^\prime$.

In the spirit of Cram\'er's model, we assume that
 an integer $n\le x$ is $k$-gleeful with probability
given by 
\begin{equation}
  \sk(x)/x.  \label{eq:skmprob}
\end{equation}

\subsection{Duplicates for $f_k(n)>1$}

Assuming (\ref{eq:skmprob}), a first try at estimating
the probability an integer $n\le x$ is a duplicate with
$f_k(n)>1$ would be simply
\begin{eqnarray*}
    \left( \frac{\sk(x)}{x} \right)^2
    & \approx &
      \frac{ k^4 x^{4/(k+1)}}{x^2 (\log x)^{2k/(k+1)}} \\
    & = &   x^{ 4/(k+1)-2} \cdot \frac{ k^4 }{ (\log x)^{2k/(k+1)}}.
\end{eqnarray*}
This probability is $o(1/x)$, unless $k<3$.
When $k=2$ we might expect the density of duplicates to be
  around $x^{1/3}/(\log x)^{4/3}$, which implies there
  are infinitely many examples.

The only potential flaw in our logic is that two different
  gleeful representations for $n$ with the same value of $k$
  must be of different lengths.
Recall that
$$
\sk(x)=\sum_{m=1}^{M(x,k)} s_{k,m}(x).
$$
Thus, a finer estimate for the probability of a duplicate
is
\begin{eqnarray*}
   && \frac{1}{x^2} \sum_{m_1=1}^{M(x,k)} s_{m_1,k}(x)
      \sum_{m_2<m_1} s_{m_2,k}(x)\\
      &=&
      \frac{1}{x^2} \sum_{m_1=1}^{M(x,k)}\sum_{m_2<m_1} s_{m_1,k}(x)
       s_{m_2,k}(x) \\
       &\approx&\frac{1}{x^2} \sum_{m_1=1}^{M(x,k)}\sum_{m_2<m_1}
       \pi( (x/m_1)^{1/k} ) \cdot  \pi( (x/m_2)^{1/k} ).
       \end{eqnarray*}
By the prime number theorem, this is asymptotic to
       \begin{eqnarray*}
       &\sim&\frac{1}{x^2} \sum_{m_1=1}^{M(x,k)}\sum_{m_2<m_1}
       \frac{k^2 x^{2/k}}{(m_1m_2)^{1/k}\log(x^2/(m_1m_2))}\\
       &\approx&
       \frac{k^2 x^{2/k-2}}{2\log(x/M)}
       \sum_{m_1=1}^{M(x,k)}m_1^{-1/k}\sum_{m_2<m_1}m_2^{-1/k}\\
       &\approx&
       \frac{k^2 x^{2/k-2}}{2\log(x/M)}
       \sum_{m_1=1}^{M(x,k)}m_1^{-1/k} \cdot\frac{m_1^{1-1/k}}{1-1/k}\\
        &\approx&
       \frac{k^3 x^{2/k-2}}{2(k-1)\log(x/M)}
       \sum_{m_1=1}^{M(x,k)}m_1^{1-2/k} .
\end{eqnarray*}
If $k=2$, the sum is just $M(x,2)\sim 3x^{1/3}/(\log x)^{2/3}$,
which gives a probability of
$$
\frac{4 M}{x\log(x/M)}
\sim
\frac{18 }{x^{2/3}(\log x)^{5/3}},
$$
smaller than our first estimate by a factor of roughly
$(\log x)^{1/3}$, but still enough that we expect infinitely
many integers $n$ with $f_2(n)>1$.

If $k\ge 3$, we end up with the following probability:
$$
\frac{k^3 x^{2/k}}{2(k-1)x^2\log(x/M)} \cdot
\frac{M^{2-2/k}}{2-2/k}.
$$
Plugging in our estimate that 
$M\approx (k+1)x^{1/(k+1)}/(\log x)^{k/(k+1)}$,
we obtain the probability
$$
\frac{k^3(k+1)}{2(k-1)} \cdot
\frac{x^{4/(k+1)-2}}{(\log x)^{(3k-1)/(k+1)}}.
$$
For $k>3$ this is clearly $o(1/x)$, as the exponent
on $x$ is less than $-1$.
For $k=3$ the exponent on $x$ is exactly $-1$, but the log factor
in the denominator still gives us $o(1/x)$.

\begin{figure}
\caption{Comparing $d(x)$ to the number of $n\le x$ with $f_2(n)>1$.}
\label{fig:dups}
    \includegraphics[width=\textwidth]{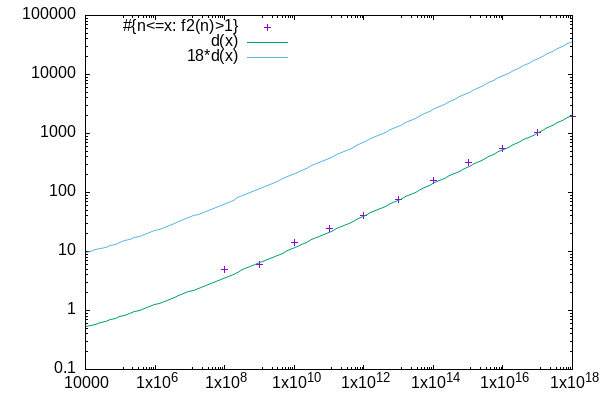}
\end{figure}

This leads us to the following conjectures.
\begin{conjecture}
    There are infinitely many integers $n$
    with $f_2(n)>1$.
\end{conjecture}
We found 1950 integers $n\le 10^{18}$ with $f_2(n)>1$.
Let us set $d(x):=x^{1/3}/(\log x)^{5/3}$, the number of integers $n$ below $x$
we expect to find with $f_2(n)>1$,
with the constant factor 18 dropped.
As you can see in Figure \ref{fig:dups}, $d(x)$ lines up very nicely with our data.
We currently have no explanation for why our prediction above is off by a factor
of $18$.

\begin{conjecture}
  For each integer $k\ge 3$,
    there are finitely many integers $n$
    with $f_k(n)>1$.
\end{conjecture}
With a bit more work, our heuristics also lead to this
much stronger conjecture:
\begin{conjecture}
    There are finitely many integers $n$
    with $f_k(n)>1$ for \textit{any} $k\ge 3$.
\end{conjecture}
We have not found any examples $n$ with $f_k(n)>1, k>2$.

Our code and data are available on the second author's github repository
at \url{https://github.com/sorenson64/sopp}.

%%%%%%%%%%%%%%%%%%%%%%%%%%%%%%%%%%%%%%%%%%%%%%%%%%%%%%%

\subsection{Duplicates for $f_k(n)>0$ and $f_{k^\prime}(n)>0$ with $k< k^\prime$}

We continue to assume $k<k^\prime$.
The heuristic probability that a randomly chosen integer $n\le x$
has both $f_k(n), f_{k^\prime}(n)>0$ is at most
$$
 \frac{s_k(x)s_{k^\prime}(x)}{x^2}
 \approx
 (kk^\prime)^2 \frac{x^{2/(k+1)+2/(k^\prime+1)-2}}{(\log x)^{2k/(k+1)+2k^\prime/(k^\prime+1)}} .
$$
With the log factors in the denominator, we can expect
infinitely many examples if the exponent on $x$ 
is strictly greater than -1, that is,
$2/(k+1)+2/(k^\prime+1)-2>-1$, or
$$\frac{2}{k+1}+\frac{2}{k^\prime+1}>1.$$
With $k^\prime>k$, this is not true when $k\ge 3$.
With $k=2$, we then require
$2/(k^\prime+1)>1/3$.
This gives $k^\prime=3$ or $4$.

\begin{conjecture}
    For $k^\prime=3$ or $4$,
    there are infinitely many integers $n$
    with both $f_2(n)>1$ and $f_{k^\prime}(n)>1$.
\end{conjecture}
We have 3 examples of 2-3 duplicates up to $10^{18}$:
\begin{quote}
    23939 \\
432958700126053 \\
137610738498311684
\end{quote}
We found no 2-4 duplicates below $10^{18}$.
We are hopeful that more examples will eventually be found.

\begin{conjecture}
   For $k<k^\prime$, 
   if $k\ge 3$ or $k^\prime\ge 5$ then
   there are finitely many integers $n$
   with both $f_k(n)>1$ and $f_{k^\prime}(n)>1$.
\end{conjecture}

\nocite{OSS2024,Moser63,TWS2022,RS62,OEIS}
\nocite{GHR,Reif}\nocite{Knuthv3}

\section{Acknowledgements}
This work was supported in part by the Mathematics Research Camp at
  Butler University in August 2023,
  and by a grant from the Holcomb Awards Committee.

Thanks to Frank Levinson for supporting Butler's research computing infrastructure.

This work began when the first author was an undergraduate student
  at Butler University.

%\bibliography{all}

\end{document}